\documentclass[12pt]{amsart}
\usepackage{amssymb}
\usepackage{amsmath}
\usepackage{amscd}
\theoremstyle{plain}
\newtheorem{theorem}{Theorem}[section]

\newtheorem{proposition}[theorem]{Proposition}
\newtheorem{lemma}[theorem]{Lemma}
{\theoremstyle{remark}

\newtheorem{remark}[theorem]{Remark}}
{\theoremstyle{definition}
\newtheorem{definition}[theorem]{Definition}
}

\newcommand{\benu}{\begin{enumerate}\renewcommand{\labelenumi}{{\rm (\roman{enumi})}}\renewcommand{\itemsep}{0pt}}
\newcommand{\eenu}{\end{enumerate}}

\newcommand{\N}{\mathbb{N}}

\newcommand{\C}{\mathbb{C}}
\newcommand{\T}{\mathbb{T}}

\newcommand{\cK}{{\mathcal K}}

\newcommand{\cH}{{\mathcal H}}

\newcommand{\K}{\mathbb{K}}

\newcommand{\fA}{\mathfrak{A}}
\newcommand{\cO}{{\mathcal O}}
\newcommand{\cT}{{\mathcal T}}

\newcommand{\s}[3]{{{#1}^{#2}_{\rm{\scriptsize #3}}}}

\newcommand{\Ca}{$C^*$-al\-ge\-bra }
\newcommand{\CA}{$C^*$-al\-ge\-bra}

\DeclareMathOperator{\id}{id}
\DeclareMathOperator{\spec}{sp}

\DeclareMathOperator{\ev}{ev}

\begin{document}
\title[$C^*$-algebras generated by scaling elements]
{\boldmath{$C^*$}-algebras generated by scaling elements}
\author[Takeshi KATSURA]{Takeshi KATSURA}
\address{Department of Mathematics, 
Hokkaido University, Kita 10, Nishi 8, 
Kita-Ku, Sapporo, 060-0810, JAPAN}
\email{katsura@math.sci.hokudai.ac.jp}



\begin{abstract}
We investigate $C^*$-algebras generated by scaling elements. 
We generalize the Wold decomposition and Coburn's theorem 
on isometries to scaling elements. 
We also completely determine 
when the $C^*$-algebra generated by a scaling element 
contains an infinite projection. 
\end{abstract}

\maketitle

\setcounter{section}{-1}

\section{Introduction}

For a (complex) Hilbert space $\cH$, 
we denote by $\mathfrak{B}(\cH)$ 
the \Ca of all bounded operators on $\cH$, 
and define a Hilbert space $\cH^\infty$ by 
$$\cH^\infty
=\bigg\{(\xi_n)_{n\in\N}\ \bigg|\ 
\xi_n\in \cH, \sum_{n\in\N}\|\xi_n\|^2<\infty\bigg\}.$$
An isometry on $\cH$ is a bounded operator $T\in \mathfrak{B}(\cH)$ 
satisfying $T^*T=I_\cH$ where $I_\cH$ is the identity operator in $\mathfrak{B}(\cH)$. 
An isometry $T$ is said to be {\em proper} if $TT^*\neq I_\cH$. 
For a Hilbert space $\cH$, 
we define an operator $S_\cH\in \mathfrak{B}(\cH^\infty)$ 
by $S_\cH(\xi_0,\xi_1,\ldots,)=(0,\xi_0,\xi_1,\ldots,)$ 
for $(\xi_0,\xi_1,\ldots,)\in \cH^\infty$. 
The operator $S_\cH$ is a proper isometry when $\cH\neq 0$. 
Conversely, 
the following proposition, called the {\em Wold decomposition}, 
says that any proper isometry is essentially 
the direct sum of an element in this form and a unitary 
(see \cite[Theorem V.2.1]{Da} for a proof). 

\medskip
\noindent
{\bfseries Proposition.} 
Let $T$ be an isometry in $\mathfrak{B}(\cH)$. 
Then there exist two Hilbert spaces $\cH_1$ and $\cH_2$, 
and an isomorphism $\cH\cong \cH_1^\infty\oplus \cH_2$ 
such that $T$ is unitarily equivalent to 
$S_{\cH_1}\oplus U\in \mathfrak{B}(\cH_1^\infty\oplus \cH_2)$ 
for a unitary $U\in \mathfrak{B}(\cH_2)$. 

The isometry $T$ is proper 
if and only if $\cH_1\neq 0$. 
\medskip

We call $S=S_\C\in \mathfrak{B}(\C^\infty)$ 
the {\em unilateral shift}. 
The $C^*$-algebra 
generated by the unilateral shift $S$ 
is called the {\em Toeplitz algebra} and denoted by $\cT$. 
In \cite{C}, 
Coburn showed the following theorem using the Wold decomposition.

\medskip
\noindent
{\bfseries Theorem.} (Coburn)
For a proper isometry $T$, 
there is a $*$-isomorphism $\varphi$ from $\cT$ 
to the \Ca $C^*(T)$ generated by $T$ such that $\varphi(S)=T$. 
\medskip

A projection $P$ in a \Ca $\fA$ is said to be {\em infinite} 
if there exists $U\in \fA$ such that $U^*U=P$ and $UU^*<P$ 
where $UU^*<P$ means $UU^*\leq P$ and $UU^*\neq P$. 
Otherwise we say that a projection $P\in\fA$ is {\em finite}. 
Existence of a proper isometry on a Hilbert space $\cH$ 
shows that the unit $I_\cH$ of $\mathfrak{B}(\cH)$ is infinite. 
It is important to determine whether a given \Ca 
contains an infinite projection or not. 
To this end, 
Blackadar and Cuntz introduced the following notion in \cite{BC}. 

\medskip
\noindent
{\bfseries Definition.} (Blackadar, Cuntz)
An element $X$ of a $C^*$-algebra is called 
a {\em scaling element} if $(X^*X)X=X$ and $X^*X\neq XX^*$. 
\medskip

Note that the condition $(X^*X)X=X$ is equivalent to 
$(X^*X)(XX^*)=XX^*$. 
Since a partial isometry $U$ satisfying $UU^*<U^*U$ 
is a scaling element, 
a \Ca containing an infinite projection 
has a scaling element. 
The converse is true if a \Ca is unital or simple 
(see \cite{BC} or \cite[Proposition 4.2]{Ka1}). 
However there exists a (non-unital, non-simple) $C^*$-algebra 
which has a scaling element but does not have an infinite projection 
(see Theorem E). 

In this paper, 
we generalize two results above on isometries 
to scaling elements, 
and investigate the structure of \CA s generated by scaling elements. 
To state the main results, 
we need several notions. 

\medskip
\noindent
{\bfseries Definition.}
For an operator $A\in \mathfrak{B}(\cH)$, 
we define an operator $S_A\in \mathfrak{B}(\cH^\infty)$ 
by $S_A(\xi_0,\xi_1,\ldots,)=(0,A\xi_0,\xi_1,\ldots,)$ 
for $(\xi_0,\xi_1,\ldots,)\in \cH^\infty$. 
$$S_A=\left(\begin{array}{ccccc}0&&&&\\
A&0&&&\\
&I_\cH&0&&\\
&&I_\cH&0&\\
&&&\hspace*{-0.3cm}\lefteqn{\ddots}&\ddots
\end{array}\right)$$
\medskip

We have $S_{I_\cH}=S_\cH$. 
For any operator $A\in \mathfrak{B}(\cH)$, 
$S_A\in \mathfrak{B}(\cH^\infty)$ is a scaling element. 
Our first result is the {\em Wold decomposition} of scaling elements. 

\medskip
\noindent
{\bfseries Theorem A.}
Let $X\in \mathfrak{B}(\cH)$ be an element satisfying $(X^*X)X=X$. 
Then there exist three Hilbert spaces $\cH_1$, $\cH_2$ and $\cH_3$, 
and an isomorphism $\cH\cong \cH_1^\infty\oplus \cH_2\oplus \cH_3$ 
such that $X$ is unitarily equivalent to 
$S_A \oplus U\oplus 0$ 
where $A$ is a positive operator in $\mathfrak{B}(\cH_1)$ 
whose support is $I_{\cH_1}$ and $U\in \mathfrak{B}(\cH_2)$ is a unitary. 
This decomposition is unique up to unitary equivalence. 

Such $X$ is a scaling element (i.e.\ $X^*X\neq XX^*$) 
if and only if $\cH_1\neq 0$. 
\medskip

We will prove Theorem A in Section 1. 
Next we state the generalization of Coburn's theomrem 
to scaling elements. 
We denote by $\spec(A)$ the spectrum of an operator $A$, 
and by $1_C$ the characteristic function of a set $C$. 
It is easy to see $0,1\in \spec(|X^*|)\subset [0,\infty)$ 
for a scaling element $X$. 

\medskip
\noindent
{\bfseries Definition.} 
A scaling element $X$ is said to be {\em non-proper} 
if $\spec(|X^*|)\setminus\{0,1\}$ is compact, 
and $1_{\{1\}}(|X|)=1_{\spec(|X^*|)\setminus \{0\}}(|X^*|)$. 
Otherwise, we say that a scaling element $X$ is {\em proper}. 
\medskip

For a positive operator $A\in \mathfrak{B}(\cH)$, 
$S_A\in \mathfrak{B}(\cH^\infty)$ is a scaling element 
with $\spec(|S_A^*|)=\spec(A)\cup\{0,1\}$. 
Using the operators $S_A$, 
we can show the following. 

\medskip
\noindent
{\bfseries Theorem B} (Existence theorem). 
For any compact set $\Omega$ with $0,1\in \Omega\subset [0,\infty)$, 
there exists a scaling element $X$ with $\spec(|X^*|)=\Omega$. 
If such $\Omega$ satisfies that 
$\Omega\setminus\{0,1\}$ is non-empty and compact, 
then there exist both a proper scaling element and a non-proper one 
whose spectra are $\Omega$. 
\medskip

Note that a scaling element $X$ is proper 
if $\spec(|X^*|)\setminus\{0,1\}$ is empty or non-compact. 
The following uniqueness theorem is a generalization of 
Coburn's theorem. 

\medskip
\noindent
{\bfseries Theorem C} (Uniqueness theorem). 
Let $X,Y$ be scaling elements.
There exists a $*$-isomorphism $\varphi\colon C^*(X)\to C^*(Y)$ 
with $\varphi(X)=Y$ if and only if $\spec(|X^*|)=\spec(|Y^*|)$ 
and $X,Y$ are simultaneously proper or non-proper. 
\medskip

Theorem C easily follows from the next proposition 
which concerns with a hierarchy of \CA s $C^*(X)$ 
generated by scaling elements $X$. 

\medskip
\noindent
{\bfseries Proposition D.}
Let $X,Y$ be scaling elements. 
\benu
\item When $X$ is proper, 
there exists a $*$-homomorphism $\varphi\colon C^*(X)\to C^*(Y)$ 
with $\varphi(X)=Y$ if and only if $\spec(|X^*|)\supset \spec(|Y^*|)$. 
\item When $X$ is non-proper, 
there exists a $*$-homomorphism $\varphi\colon C^*(X)\to C^*(Y)$ 
with $\varphi(X)=Y$ if and only if $\spec(|X^*|)\supset \spec(|Y^*|)$ 
and $Y$ is non-proper. 
\eenu
\medskip

Clearly the $*$-homomorphism $\varphi$ in the proposition above is, 
if it exists, unique and surjective. 
In Section 3 
we will prove Proposition D and the following theorem, 
by using the theory of \CA s arising from constant maps 
studied in Section 2. 

\medskip
\noindent
{\bfseries Theorem E.}
Let $X$ be a scaling element.
The \Ca $C^*(X)$ has an infinite projection 
if and only if $[0,1]\setminus\spec(|X^*|)\neq\emptyset$. 
\medskip

\section{The Wold decomposition of scaling elements}

In this section, 
we prove the Wold decomposition of scaling elements (Theorem A). 
For an operator $X\in \mathfrak{B}(\cH)$, we denote by $l(X)$ (resp.\ $r(X)$) 
the left (resp. right) support of $X$. 
Namely $l(X)$ (resp.\ $r(X)$) is the smallest projection of $\mathfrak{B}(\cH)$ 
satisfying $l(X)X=X$ (resp.\ $Xr(X)=X$). 

Take an element $X\in \mathfrak{B}(\cH)$ satisfying $(X^*X)X=X$. 
We set $P_0=r(X)$ and $P_0'=l(X)$. 
Since $(X^*X)X=X$, 
we have $X^*XP_0'=P_0'$. 
Hence we get $P_0\geq P_0'$. 
Set $P_3=I_\cH-P_0$ and $\cH_3=P_3\cH$. 
Then we have $P_3X=XP_3=0$. 
Set $Q_0=P_0-P_0'$ and $X_0=XQ_0$. 
Since $r(X)=P_0\geq Q_0$, we have $r(X_0)=Q_0$. 
We define $Q_1=l(X_0)$ and $U_1=XQ_1$. 
Clearly we have $Q_1\leq P_0'$. 
Since $X^*XP_0'=P_0'$, we have $X^*XQ_1=Q_1$. 
Thus $U_1^*U_1=Q_1$. 
Since $X_0^*XP_0'=Q_0X^*XP_0'=Q_0P_0'=0$, 
we have $Q_1XP_0'=0$. 
This implies that 
$$Q_1X=Q_1XP_0=Q_1XQ_0+Q_1XP_0'=Q_1XQ_0=XQ_0.$$ 
Recursively we define projections $Q_2,Q_3,\ldots$ and 
partial isometries $U_2,U_3,\ldots$ 
by $Q_n=U_{n-1}U_{n-1}^*$ and $U_n=XQ_n$. 
From $Q_n\leq P_0'$, 
we have $X^*XQ_n=Q_n$. 
Hence $U_n^*U_n=Q_n$. 
By the definition, we have $Q_{n+1}=XQ_nX^*$ 
for $n\geq 1$. 
Thus we get $Q_{n+1}X=XQ_nX^*X=XQ_n$. 
For $n\geq 1$, 
we have $Q_0Q_n=0$ because $Q_n\leq P_0'$. 
For $k\geq 1$, 
we have $Q_kQ_{n+k}=(XQ_{k-1}X^*)(XQ_{n+k-1}X^*)=XQ_{k-1}Q_{n+k-1}X^*$. 
Thus recursively 
we can show that $Q_mQ_n=0$ for $0\leq m<n$. 
Therefore $\{Q_n\}_{n\in\N}$ is a family of mutually orthogonal projections. 
Set $P_1=\sum_{n\in\N}Q_n$. 
We have 
$$P_1X=\sum_{n=0}^\infty Q_nX=\sum_{n=1}^\infty Q_nX
=\sum_{n=1}^\infty XQ_{n-1}=XP_1.$$ 
Let $X_0=U_0|X_0|$ be a polar decomposition of $X_0$.
Then we have $U_0^*U_0=r(X_0)=Q_0$ and $U_0U_0^*=l(X_0)=Q_1$. 
Since $r(|X_0|)=l(|X_0|)=Q_0$, 
the restriction of $|X_0|$ 
on the Hilbert space $\cH_1=Q_0\cH$ 
gives a positive operator $A\in \mathfrak{B}(\cH_1)$ with $r(A)=l(A)=I_{\cH_1}$. 
By using partial isometries $\{U_n\}_{n\in\N}$, 
we have a unitary from $P_1\cH$ to $\cH_1^\infty$. 
Via this unitary, 
the operator $P_1X=XP_1\in \mathfrak{B}(P_1\cH)$ 
is unitarily equivalent to $S_A\in \mathfrak{B}(\cH_1^\infty)$. 
We set $P_2=I_\cH-P_1-P_3$ and $X_2=P_2X=XP_2$. 
Since $P_2\leq P_0'\leq P_0$, 
we have $l(X_2)=r(X_2)=P_2$. 
We can easily check $(X_2^*X_2)X_2=X_2$. 
This shows that $X_2^*X_2=X_2^*X_2P_2=P_2$. 
Hence we get $X_2^*X_2=X_2X_2^*=P_2$. 
Therefore the restriction of $X$ on $\cH_2=P_2\cH$ 
is a unitary in $\mathfrak{B}(\cH_2)$. 
The uniqueness of this decomposition follows 
from the argument above and the uniqueness of polar decomposition. 
This shows the former part of Theorem A, 
and the latter part is obvious. 

\medskip

Similarly as \cite{C}, 
we can deduce Proposition D from Theorem A. 
However the computations we need here is much harder 
than the ones in \cite{C}. 
In this paper, 
we will prove Proposition D 
by using the general theory of \CA s 
arising from constant maps.

\section{$C^*$-algebras arising from constant maps}

In this section, 
we study the structures of \CA s arising from constant maps. 
In the next section, 
we apply the results of this section 
to the \CA s generated by scaling elements. 

Take a locally compact space $\Omega$ 
and a point $v$ of $\Omega$. 
We fix them throughout this section. 
An {\em $(\Omega,v)$-pair} is a pair $(\pi,t)$ consisting 
of a $*$-homomorphism $\pi$ and a linear map $t$ 
from $C_0(\Omega)$ to some \Ca $\fA$ satisfying 
$$t(f)^*t(g)=\pi(\overline{f}{g}),\qquad
t(f)\pi(g)=t(fg),\qquad
\pi(f)t(g)=f(v)t(g)$$
for $f,g\in C_0(\Omega)$. 
Note that the second condition automatically follows 
from the first one (see \cite{Ka2}). 

\begin{definition}
Let us denote by $(\hat{\pi},\hat{t})$ 
the universal $(\Omega,v)$-pair, 
and by $\cT(\Omega,v)$ the \Ca 
generated by the images of $\hat{\pi}$ and $\hat{t}$. 
\end{definition}

The universality means that for any $(\Omega,v)$-pair $(\pi,t)$ 
on a \Ca $\fA$, there exists a $*$-homomorphism 
$\rho\colon \cT(\Omega,v)\to \fA$ 
with $\pi=\rho\circ \hat{\pi}$ and $t=\rho\circ \hat{t}$. 
The existence of the universal $(\Omega,v)$-pair follows 
from a standard argument, 
and the uniqueness is easy to see. 

\begin{lemma}\label{repre}
Both $\hat{\pi}$ and $\hat{t}$ are injective. 
\end{lemma}

\begin{proof}
It suffices to find one $(\Omega,v)$-pair $(\pi,t)$ 
such that $\pi$ and $t$ are injective. 
Let us take a faithful representation $\varphi\colon C_0(\Omega)\to \mathfrak{B}(\cH)$. 
We define an injective $*$-homomorphism 
$\pi\colon C_0(\Omega)\to \mathfrak{B}(\cH^\infty)$ by 
$$\pi(f)(\xi_0,\xi_1,\xi_2,\ldots,)
=(\varphi(f)\xi_0,f(v)\xi_1,f(v)\xi_2,\ldots,)$$ 
for $f\in C_0(\Omega)$ and $(\xi_0,\xi_1,\ldots,)\in \cH^\infty$. 
We define an injective linear map 
$t\colon C_0(\Omega)\to \mathfrak{B}(\cH^\infty)$ 
by $t(f)=S_\cH\pi(f)$ for $f\in C_0(\Omega)$. 
Then the pair $(\pi,t)$ is an $(\Omega,v)$-pair. 
We are done. 
\end{proof}

For a while, 
we suppose that $\Omega\setminus\{v\}$ is compact. 
We denote by $V$ an isometry $\hat{t}(1_{\Omega})\in\cT(\Omega,v)$. 
The \Ca $\cT(\Omega,v)$ is generated by the isometry $V$ 
and $\hat{\pi}(C(\Omega))$. 
By the proof of Lemma \ref{repre}, 
the projection $P=\hat{\pi}(1_{\{v\}})-VV^*$ 
is non-zero. 
Let us define $I\subset \cT(\Omega,v)$ 
by the closure of the linear span of 
$${\mathcal E}=\{V^nP(V^*)^m\mid n,m\in\N\}.$$

\begin{lemma}\label{I=I(P)}
The subset $I$ is the ideal of $\cT(\Omega,v)$ generated by $P$, 
\end{lemma}

\begin{proof}
It suffices to see that the linear span of ${\mathcal E}$ 
is closed under the multiplication by $V$, $V^*$ 
and $\hat{\pi}(C(\Omega))$ from left. 
For $f\in C(\Omega)$, 
we have $\hat{\pi}(f)V=f(v)V$ and 
$\hat{\pi}(f)P=\hat{\pi}(f1_{\{v\}})-\hat{\pi}(f)VV^*=f(v)P$. 
Hence the linear span of ${\mathcal E}$ 
is closed under the multiplication by $\hat{\pi}(C(\Omega))$ from left. 
Clearly it is closed under the multiplication by $V$ from left. 
We have 
$V^*P=V^*\hat{\pi}(1_{\{v\}})-V^*VV^*=V^*-V^*=0$. 
This shows that the linear span of ${\mathcal E}$ 
is also closed under the multiplication by $V^*$ from left. 
We are done. 
\end{proof}

\begin{definition}
When $\Omega\setminus\{v\}$ is compact, 
we define a \Ca $\cO(\Omega,v)$ by $\cO(\Omega,v)=\cT(\Omega,v)/I$. 
\end{definition}

We denote by $\cK$ the \Ca of compact operators 
on the Hilbert space $\C^\infty$. 
The matrix unit of $\cK$ is denoted by 
$\{U_{n,m}\mid n,m\in\N\}$. 

\begin{lemma}\label{I=K}
The map $V^nP(V^*)^m\mapsto U_{n,m}$ 
induces an isomorphism $I\cong\cK$. 
\end{lemma}

\begin{proof}
This follows from routine computation. 
\end{proof}

Let $(\pi,t)$ be an $(\Omega,v)$-pair, 
and $\fA$ be the \Ca generated by the images of $\pi$ and $t$. 
By the universality, 
there exists a surjection $\rho\colon \cT(\Omega,v)\to \fA$ 
with $\pi=\rho\circ \hat{\pi}$ and $t=\rho\circ \hat{t}$. 
When $\Omega\setminus\{v\}$ is compact 
and $\pi(1_{\{v\}})=t(1_{\Omega})t(1_{\Omega})^*$, 
the surjection $\rho\colon \cT(\Omega,v)\to \fA$ factors through 
a surjection $\bar{\rho}\colon \cO(\Omega,v)\to \fA$. 
By using results in \cite{Ka2}, 
we get the following. 

\begin{proposition}\label{univ}
\benu
\item When $\Omega\setminus\{v\}$ is not compact, 
the surjection $\rho$ is an isomorphism if and only if 
$\pi$ is injective. 
\item When $\Omega\setminus\{v\}$ is compact, 
the surjection $\rho$ is an isomorphism if and only if 
$\pi$ is injective and 
$\pi(1_{\{v\}})\neq t(1_{\Omega})t(1_{\Omega})^*$. 
\item When $\Omega\setminus\{v\}$ is a non-empty compact set 
and $\pi(1_{\{v\}})=t(1_{\Omega})t(1_{\Omega})^*$, 
the surjection $\bar{\rho}\colon \cO(\Omega,v)\to \fA$ 
is an isomorphism if and only if 
$\pi$ is injective. 
\eenu
\end{proposition}

\begin{proof}
We define a continuous map $r\colon\Omega\to \Omega$ 
by $r(x)=v$ for all $x\in \Omega$, 
and set a topological graph $E=(\Omega,\Omega,\id,r)$ 
(see \cite{Ka2}). 
If $\Omega\setminus\{v\}$ is not compact, 
then we have $\s{E}{0}{rg}=\emptyset$. 
Thus $\cT(\Omega,v)=\cT(E)=\cO(E)$. 
When $\Omega\setminus\{v\}$ is compact, 
we have $\s{E}{0}{rg}=\{v\}$. 
Hence we get $\cT(\Omega,v)=\cT(E)$ and $\cO(\Omega,v)=\cO(E)$. 
It is easy to verify that 
the topological graph $E$ is topologically free 
when $\Omega\neq\{v\}$. 
Therefore (i) and (iii) follows from \cite[Theorem 5.12]{Ka2}, 
and (ii) follows from \cite[Proposition 3.16]{Ka22}. 
\end{proof}

When $\Omega=\{v\}$, 
we have $\cT(\Omega,v)\cong\cT$ and $\cO(\Omega,v)\cong C(\T)$ 
where $\T$ is the one-dimensional torus. 
Thus in this case, (ii) in the proposition above 
is Coburn's Theorem introduced in the introduction, 
and the corresponding statement of (iii) does not hold. 

Let us take a closed subset $\Omega'$ of $\Omega$ with $v\in \Omega'$. 
The universal $(\Omega',v)$-pair is denoted by $(\hat{\pi}',\hat{t}')$. 
By the universality, 
there exists a surjection $\rho\colon \cT(\Omega,v)\to \cT(\Omega',v)$ 
with $\hat{\pi}'\circ\sigma=\rho\circ\hat{\pi}$ 
and $\hat{t}'\circ\sigma=\rho\circ \hat{t}$ 
where $\sigma\colon C_0(\Omega)\to C_0(\Omega')$ is 
the natural surjection. 
The kernel of $\sigma$ is $C_0(\Omega\setminus \Omega')$. 

\begin{proposition}
The kernel of $\rho$ is isomorphic to $C_0(\Omega\setminus \Omega')\otimes\cK$. 
\end{proposition}

\begin{proof}
Let $J$ be the closure of the linear span of elements in the form 
$$\hat{t}(f_1)\hat{t}(f_2)\cdots \hat{t}(f_n)\hat{\pi}(h)\hat{t}(g_m)^*\cdots \hat{t}(g_2)^*\hat{t}(g_1)^*$$
where 
$f_1,\ldots,f_n,g_1,\ldots,g_m\in C_0(\Omega)$ 
and $h\in C_0(\Omega\setminus \Omega')$. 
It is routine to check that $J$ is the ideal generated 
by $\hat{\pi}(C_0(\Omega\setminus \Omega'))$. 
It is also routine to see that the map 
\begin{align*}
J\ni \hat{t}(f_1)&\cdots \hat{t}(f_n)\hat{\pi}(h)
\hat{t}(g_m)^*\cdots \hat{t}(g_1)^*\\
\mapsto& f_1(v)\cdots f_n(v)
\overline{g_m(v)}\cdots \overline{g_1(v)}
\big(h\otimes U_{n,m}\big)\in C_0(\Omega\setminus \Omega')\otimes\cK
\end{align*}
induces an isomorphism. 
We will show that $J=\ker\rho$. 
Since $\hat{\pi}(C_0(\Omega\setminus \Omega'))\subset \ker\rho$, 
we have $J\subset \ker\rho$. 
Hence the surjection $\rho\colon \cT(\Omega,v)\to \cT(\Omega',v)$ 
factors through $\varphi\colon \cT(\Omega,v)/J\to \cT(\Omega',v)$. 
Since 
$\hat{\pi}(C_0(\Omega\setminus \Omega')), 
\hat{t}(C_0(\Omega\setminus \Omega'))\subset J$, 
we have a $*$-homomorphism 
$\pi'\colon C_0(\Omega')\to \cT(\Omega,v)/J$ 
and a linear map $t'\colon C_0(\Omega')\to \cT(\Omega,v)/J$ 
such that $\hat{\pi}'=\varphi\circ\pi'$ and $\hat{t}'=\varphi\circ t'$. 
The pair $(\pi',t')$ is an $(\Omega',v)$-pair. 
Hence we get a $*$-homomorphism $\cT(\Omega',v)\to \cT(\Omega,v)/J$ 
which is clearly the inverse of the surjection $\varphi$. 
Thus $J=\ker(\rho)$. 
\end{proof}

When $\Omega\setminus\{v\}$ is compact, 
$\Omega'\setminus\{v\}$ is also compact. 
In this case, 
the surjection $\rho\colon \cT(\Omega,v)\to \cT(\Omega',v)$ 
induces the $*$-homomorphism 
$\bar{\rho}\colon \cO(\Omega,v)\to \cO(\Omega',v)$. 
By Lemma \ref{I=I(P)}, the kernel of the surjection $\bar{\rho}$
is also isomorphic to $C_0(\Omega\setminus \Omega')\otimes\cK$. 

By taking $\Omega'=\{v\}$, 
we get the following commutative diagram 
with two exact rows;
$$\begin{CD}
0 @>>> C_0(\Omega\setminus \{v\}) @>>> C_0(\Omega) @>\ev_v >>\C @>>> 0 \\
@. @VVV @VV\hat{\pi} V  @VVV \\
0 @>>> C_0(\Omega\setminus \{v\})\otimes\K 
@>>> \cT(\Omega,v) @>\rho >> \cT @>>> 0\lefteqn{.} \\
\end{CD}$$
By this diagram, 
we see that $\cT(\Omega,v)$ is a type I \CA . 
By Lemma \ref{I=K}, $\cO(\Omega,v)$ is also type I 
when $\Omega\setminus \{v\}$ is compact. 

\begin{proposition}\label{KKisom}
The $*$-homomorphism $\hat{\pi}\colon C_0(\Omega)\to \cT(\Omega,v)$ 
induces an isomorphism between $K$-groups. 
\end{proposition}

\begin{proof}
It is well-known that the left and right vertical maps 
in the diagram above induce isomorphisms on $K$-groups. 
Hence the 5-lemma shows that $\hat{\pi}\colon C_0(\Omega)\to \cT(\Omega,v)$ 
also induces an isomorphism between $K$-groups. 
\end{proof}

\begin{remark}
By using the isomorphism $\cT(\Omega,v)=\cT(E)$ 
explained in the proof of Proposition \ref{univ}, 
Proposition \ref{KKisom} follows from \cite[Lemma 6.5]{Ka2} 
(see also \cite[Theorem 4.4]{P} or \cite[Proposition 8.2]{Ka3}). 
By \cite[Corollary 6.10]{Ka2} with some computation, 
we see that the $K$-groups of the \Ca $\cO(\Omega,v)=\cO(E)$ 
is isomorphic to the ones of $C(\Omega\setminus\{v\})$ 
when $\Omega\setminus\{v\}$ is compact and non-empty. 
\end{remark}

The next two lemmas are standard, hence we omit the proofs. 

\begin{lemma}\label{LemK1}
Let $\Omega$ be a compact space, 
and $v,w$ be points in $\Omega$. 
The two evaluation maps $\ev_v,\ev_w\colon C(\Omega)\to\C$ 
induce same maps on $K$-groups 
if and only if $v$ and $w$ are in a same connected component. 
\end{lemma}

\begin{lemma}\label{LemK2}
A projection $P$ of $\cT$ defines $0$ in $K_0(\cT)$ 
if and only if $P$ is finite. 
\end{lemma}

\begin{proposition}\label{infinite}
The \Ca $\cT(\Omega,v)$ has an infinite projection 
if and only if there exists a compact open subset $C$ of $\Omega$ 
containing $v$. 
\end{proposition}

\begin{proof}
Suppose that there exists a compact open subset $C$ of $\Omega$ 
containing $v$. 
Set $U=\hat{t}(1_{C})\in \cT(\Omega,v)$.  
We will show that this partial isometry $U$ 
satisfies $UU^*<U^*U$. 
When $\Omega=\{v\}$, the \Ca $\cT(\Omega,v)$ 
is isomorphic to the Toeplitz algebra 
and $U$ is a proper isometry. Thus $UU^*<U^*U$.
Otherwise, we can find $w_1\in C\setminus\{v\}$ 
or $w_2\in \Omega\setminus C$. 
When there exists $w_1\in C\setminus\{v\}$, 
we can find $f\in C_0(\Omega)$ satisfying $0\leq f\leq 1$, $f(v)=1$, $f(w_1)=0$ 
and $f(w)=0$ for $w\in \Omega\setminus C$. 
Then we have 
$$UU^*=\hat{\pi}(f)UU^*\hat{\pi}(f)\leq \hat{\pi}(f^2)<\hat{\pi}(1_{C})=U^*U.$$
Hence $UU^*<U^*U$.
When there exists $w_2\in \Omega\setminus C$, 
we can find $g\in C_0(\Omega)$ satisfying $0\leq g\leq 1$, $g(w_2)=1$, 
and $g(w)=1$ for $w\in C$. 
Then we have 
$UU^*=\hat{t}(g)\hat{\pi}(1_{C})\hat{t}(g)^*\leq \hat{t}(g)\hat{t}(g)^*$. 
Take $h\in C_0(\Omega)$ with $h(w_2)=1$ and $h(w)=0$ for $w\in C$. 
Since $UU^*\hat{t}(h)=0$ and $\hat{t}(g)\hat{t}(g)^*\hat{t}(h)=\hat{t}(g\overline{g}h)\neq 0$, 
we have $UU^*\neq \hat{t}(g)\hat{t}(g)^*$. 
Hence $UU^*<\hat{t}(g)\hat{t}(g)^*\leq U^*U$. 
Therefore, if there exists a compact open subset $C$ of $\Omega$ 
containing $v$ 
then the \Ca $\cT(\Omega,v)$ has an infinite projection. 

Conversely suppose that there exist no compact open subsets of $\Omega$ 
containing $v$. 
To derive a contradiction, 
we assume that the \Ca $\cT(\Omega,v)$ 
has an infinite projection $P$. 
Take $U\in \cT(\Omega,v)$ with $U^*U=P$ and $UU^*<P$. 
Set $Q=UU^*$ and $P_0=P-Q$. 
By Lemma \ref{LemK1}, 
the evaluation map $\ev_v\colon C_0(\Omega)\to\C$ at $v\in \Omega$ 
induces $0$ map between $K$-groups. 
Hence the natural surjection $\rho\colon\cT(\Omega,v)\to\cT$ also 
induces $0$ map between $K$-groups, 
because the three vertical maps in the following commutative diagram 
induce isomorphisms on $K$-groups; 
$$\begin{CD}
0 @>>> C_0(\Omega\setminus \{v\}) @>>> C_0(\Omega) @>\ev_v >>\C @>>> 0 \\
@. @VV\hat{\pi} V @VV\hat{\pi} V  @VVV \\
0 @>>> \ker\rho 
@>>> \cT(\Omega,v) @>\rho >> \cT @>>> 0\lefteqn{.} \\
\end{CD}$$
Thus $\rho(P)$ defines $0$ in $K_0(\cT)$. 
By Lemma \ref{LemK2}, $\rho(P)$ is finite. 
Hence we have $\rho(P)=\rho(Q)$. 
Thus $P_0\in \ker\rho$. 
Since the map $\hat{\pi}\colon C_0(\Omega\setminus \{v\})\to \ker\rho$ is 
an isomorphism onto a full corner, 
there exists a non-zero projection $P_0'$ in $M_n(C_0(\Omega\setminus \{v\}))$ 
for some $n\in\N$ 
whose image by $\hat{\pi}$ defines the same element as $P_0$ 
in $K_0(\ker\rho)$. 
Since every non-zero projection in $M_n(C_0(\Omega))$ 
defines non-zero elements in $K_0(C_0(\Omega))$, 
the image of $P_0'$ 
by the natural embedding $C_0(\Omega\setminus \{v\})\to C_0(\Omega)$ 
defines a non-zero element in $K_0(C_0(\Omega))$. 
This shows that $P_0\in \cT(\Omega,v)$ 
defines a non-zero element in $K_0(\cT(\Omega,v))$. 
This is a contradiction because $P_0=U^*U-UU^*$. 
Therefore if there exist no compact open subsets of $\Omega$ 
containing $v$, then $\cT(\Omega,v)$ has no infinite projections. 
\end{proof}

\begin{remark}
The proof of Proposition \ref{infinite} 
shows that the \Ca $\cT(\Omega,v)$ is stably finite 
when there exist no compact open subsets of $\Omega$ 
containing $v$. 
However its unitization has an infinite projection 
because it has a scaling element. 
Hence the \Ca $\cT(\Omega,v)$ is not quasi-diagonal. 
\end{remark}

\begin{proposition}\label{infinite2}
When $\Omega\setminus\{v\}$ is a non-empty compact set, 
the \Ca $\cO(\Omega,v)$ has an infinite projection. 
\end{proposition}

\begin{proof}
Let $U$ be the image of $\hat{t}(1_{\{v\}})\in \cT(\Omega,v)$ 
via the natural surjection $\cT(\Omega,v)\to \cO(\Omega,v)$. 
Then similarly as in the proof of Proposition \ref{infinite}, 
we can show that $U$ is a partial isometry 
with $UU^*<U^*U$. 
\end{proof}

When $\Omega=\{v\}$, the \Ca $\cO(\Omega,v)\cong C(\T)$ 
has no infinite projections.

\section{$C^*$-algebras generated by scaling elements}

In the last section, 
we will show Proposition D and Theorem E 
by using results in the previous section. 

Take a scaling element $X$ in some \Ca $\fA$. 
Set $\Omega=\spec(|X^*|)\setminus\{0\}$ which is a locally compact space. 
We define an injective $*$-isomorphism $\pi\colon C_0(\Omega)\to \fA$ by 
$\pi(f)=f(|X|)$ for $f\in C_0(\Omega)$. 
Take $g\in C_0(\Omega)$ 
and define $f\in C_0(\Omega)$ by $f(x)=xg(x)$. 
Then the formula $t(f)=Xg(|X|)$ extends 
a well-defined linear map $t\colon C_0(\Omega)\to \fA$ 
satisfying $t(f)^*t(g)=\pi(\overline{f}{g})$ and 
$t(f)\pi(g)=t(fg)$ for $f,g\in C_0(\Omega)$. 
Since $X$ satisfies $(X^*X)X=X$, 
we have $\pi(f)t(g)=f(1)t(g)$ for $f,g\in C_0(\Omega)$. 
Hence the pair $(\pi,t)$ is an $(\Omega,1)$-pair. 
When $\Omega$ is compact, 
we have $t(1_{\Omega})t(1_{\Omega})^*=1_{\Omega}(|X^*|)$. 
Therefore Proposition \ref{univ} shows that 
the \Ca $C^*(X)$ generated by the scaling element $X$ 
is isomorphic to $\cT(\Omega,1)$ if $X$ is proper, 
and to $\cO(\Omega,1)$ if $X$ is non-proper. 
This completes the proof of Proposition D. 
Theorem E follows from Proposition \ref{infinite} and 
Proposition \ref{infinite2}. 
We also see that the \Ca $C^*(X)$ is type I 
and that $\pi\colon C_0(\Omega)\to C^*(X)$ 
induces an isomorphism between $K$-groups 
if $X$ is proper. 
If a scaling element $X$ is non-proper, 
the $K$-groups of the \Ca $C^*(X)$ are isomorphic to 
the ones of $C(\Omega\setminus\{1\})$. 

\medskip

{\bf Acknowledgments.} 
A part of this work was done 
while the author was staying at the University of Nevada, Reno. 
He would like to thank people there for their warm hospitality. 
He is also grateful 
to Yasuyuki Kawahigashi for their encouragement. 
This work was partially supported by Research Fellowship 
for Young Scientists of the Japan Society for the Promotion of Science.

\end{document}